%% file: main.tex
\documentclass[journal,twoside,web]{ieeecolor}
\usepackage{generic}
\usepackage{cite}
\usepackage[commentmarkup=footnote]{changes}
\usepackage{amsmath,amssymb,amsfonts}
\usepackage{mathtools, booktabs}
\usepackage{algorithmic}
\usepackage{graphicx}
\usepackage{algorithm,algorithmic}
\usepackage{hyperref}
\hypersetup{hidelinks=true}
\usepackage{textcomp}
\usepackage{hyperref}

\usepackage{pgfplots}
\pgfplotsset{compat=newest}

\newcommand{\term}[1]{\emph{#1}}
\definecolor{TolVibrantBlue}{HTML}{0077BB}
\definecolor{TolVibrantCyan}{HTML}{33BBEE}
\definecolor{TolVibrantTeal}{HTML}{009988}
\definecolor{TolVibrantOrange}{HTML}{EE7733}
\definecolor{TolVibrantRed}{HTML}{CC3311}
\definecolor{TolVibrantMagenta}{HTML}{EE3377}
\definecolor{TolVibrantGray}{HTML}{BBBBBB}

\definecolor{TolMutedBlue}{HTML}{332288}
\definecolor{TolMutedCyan}{HTML}{88CCEE}
\definecolor{TolMutedTeal}{HTML}{44AA99}
\definecolor{TolMutedGreen}{HTML}{117733}
\definecolor{TolMutedOlive}{HTML}{999933}
\definecolor{TolMutedSand}{HTML}{DDCC77}
\definecolor{TolMutedRose}{HTML}{CC6677}
\definecolor{TolMutedWine}{HTML}{882255}
\definecolor{TolMutedPurple}{HTML}{AA4499}

\definechangesauthor[name = {Ronny Bergmann}, color = TolVibrantMagenta]{RB}
\definechangesauthor[name = {Mateusz Baran}, color = TolVibrantOrange]{MB}

\NewDocumentCommand{\deriv}{ m m }{\frac{\mathrm{d}}{\mathrm{d}#1}{#2}}

\def\BibTeX{{\rm B\kern-.05em{\sc i\kern-.025em b}\kern-.08em
    T\kern-.1667em\lower.7ex\hbox{E}\kern-.125emX}}
\begin{document}

\title{A geometric perspective of state estimation using Kalman filters}
\author{Mateusz Baran, Ronny Bergmann
\thanks{Mateusz Baran is with AGH University of Krakow,
30 Mickiewicz Ave., 30-059 Krakow, Poland (e-mail: mbaran@agh.edu.pl). }
\thanks{Ronny Bergmann is with Norwegian University of Science and Technology, Department of Mathematical Sciences, NO-7041 Trondheim, Norway
(e-mail: ronny.bergmann@ntnu.no).}
}

\maketitle

\begin{abstract}
Geometry of the state space is known to play a crucial role in many applications of Kalman filters, especially robotics and motion tracking.
The Lie group-centric approach is currently very common, although a Riemannian approach has also been developed.
In this work we explore the relationship between these two approaches and develop a novel description of Kalman filters based on affine connections that generalizes both commonly encountered descriptions.
We illustrate the results on two test problems involving the special Euclidean group and the tangent bundle of a sphere in which the state is tracked by geometric variants of the extended Kalman filter and the unscented Kalman filter.
The examples use a newly developed library GeometricKalman.jl.
The new approach provides a greater freedom in selecting the structure of the state space for state estimation and can be easily integrated with standard techniques such as parameter estimation or covariance matrix estimation.
\end{abstract}

\begin{IEEEkeywords}
computational geometry, filtering theory, Kalman filters, Lie groups, nonlinear filters
\end{IEEEkeywords}

\input{content.tex}


\section*{References}

\bibliographystyle{IEEEtran}
\bibliography{manuscript}

\end{document}

%% file: content.tex
\section{Introduction}
\label{sec:introduction}

Over the decades numerous different variants of the Kalman filter have been developed~\cite{BarrauBonnabel:2018,MarinsXiaopingBachmannMcGheeZyda:2001,Im:2024,RigoSegalaSansonettoMuradore:2022}.
Many of them find applications in robotics or navigation, where the geometry of the state space is known to play a crucial role~\cite{Groves:2013}.

The core subject of our study are dynamical systems affected by noise.
We would like to estimate parameters of state of such systems based on our knowledge of its dynamics and some measurements.
In the continuous and deterministic setting it can be formalized as tracking the state of a system $p(t)$ which evolves in the configuration space $\mathcal M$.
The initial state is $p_0 = p(0) \in \mathcal M$.
This evolution is described by a differential equation
\begin{equation}
    \label{eq:f_system}
    \deriv{t}{p(t)} = f(p(t), q(t), t)
\end{equation}
for some function $f$.
Its first argument, $p(t)$, is the state at time $t \in [0, T]$, $q(t)$ describe control parameters and $T$ is the upper bound on simulation time.
An example of such system can be found in~\cite{BrossardBarrauChauchatBonnabel:2021}.
Kalman filters allow us to estimate the value of $p$ in situations where there is significant noise, $f$ might be a very rough approximation to the true dynamics and we have some observations to help us track the state.

In this work we describe a novel and generalized presentation of Kalman filters.
In this framework we explain how similarities of multiple existing approaches like error-state Kalman filters, invariant filter and optimization-based filters can be emphasized and used to develop new variants.
We explain theoretical foundations in a way that should be approachable by practitioners without omitting important parts of the theory.
We attempted to collect the most prominent variants of Kalman filters and present them in a unified geometric framework.
Lie groups are commonly used in robotics and aeronautics~\cite{BarrauBonnabel:2018} as well as in control theory~\cite{BulloLewis:2005} to represent configuration spaces.
A recent literature review can be found in~\cite{MahonyGoorHamel:2022}.
Riemannian manifolds have been used for example for articulated tracking~\cite{HaubergSommerPedersen:2012} and more generally in computer vision~\cite{TuragaSrivastava:2016}.
We generalize both approaches to arbitrary manifolds with an affine connection.
We demonstrate our approach by generic implementation of selected variants \verb|GeometricKalman.jl| (available at \url{https://github.com/JuliaManifolds/GeometricKalman.jl}) based on one of the Julia packages for numerical differential geometry, \verb|Manifolds.jl|~\cite{AxenBaranBergmannRzecki:2023} and Riemannian optimization, \verb|Manopt.jl|~\cite{Bergmann:2022}.

This work was enabled by recent developments in applications of differential geometry, such as computational anatomy~\cite{PennecLorenzi:2020,kuhnel_computational_2017}, shape analysis~\cite{srivastava_shape_2011,baran_closest_2018,Rzecki:2021} or medical image analysis~\cite{BergmannFitschenPerschSteidl:2018:1,BacakBergmannSteidlWeinmann:2016:1}.
Other notable achievements include the Geomstats library~\cite{miolane_geomstats_2020}, which implements selected variants of Kalman filters on manifolds.

As a guiding example we consider a segmented rigid body.
Its state is described by the position and orientation of the first segment as well as pairs of angles between two consecutive segments.
Such model is widely considered in robotics~\cite{HeWangLiu:2019} and modelling of human and animal skeletons, especially the hip and shoulder joints~\cite{Roupa:2022}.
The natural geometry of those pairs of angles is that of a sphere~\cite{CelledoniCokajLeoneMurariOwren:2021}, which cannot be handled using the standard approach of Lie group representation~\cite{FornasierGoorAllakMahonyWeiss:2023}.

Following contributions are made in this work.
First, we demonstrate that the theory of connection manifolds is more fundamental to Kalman filters than Lie theory.
Next, we provide a systematic review of different variants of Kalman filters in the developed framework.
A new library using this approach was written.
Finally, we also show some paths for future developments this approach opens.

\newpage
\section{Theoretical foundations}

\subsection{Basic concepts in differential geometry and Lie groups}

Manifolds with an affine connection~\cite{PennecLorenzi:2020,Lee:2012:1} are a commonly encountered basic differential geometry objects useful in applications.
An $d$-dimensional \term{connection manifold} $\mathcal{M}$ is characterized by an atlas of charts $\varphi_i \colon U_i \to \mathbb{R}^d$, $i \in \mathcal I$, where $\mathcal I$ is some index set, $U_i \subseteq \mathcal{M}$ is open, and $\varphi_i$ is a smooth bijection. The union of all $U_i$, $i\in\mathcal I$ is $\mathcal M$.
At each point we have a \term{tangent space} $T_p \mathcal{M}$, which is a vector space defined as follows:
In a chart $(\varphi, U)$ around $p \in U$ consider curve $\gamma\colon I \to U\subset\mathcal M$, where $I$ an interval containing $0$, such that $\gamma(0) = p$.
We call the equivalence class $X = [\deriv{t}{\varphi\circ\gamma}] \in T_p\mathcal M$ of curves $\delta$ through $\delta(0)$ having the same derivative with respect to $\varphi$ as $\gamma$ a \term{tangent vector}. This definition is indeed independent of the chart $\varphi$, see also~\cite[p.~33]{AbsilMahonySepulchre:2008}.

We represent $X$ also in coordinates with respect to the chart $\varphi$, by writing the curve $\gamma$ as
\begin{equation*}
    \gamma(t) = \varphi^{-1}\bigl(
        \varphi(p) + tc
    \bigr),\qquad \text{ for some } c \in \mathbb R^d
\end{equation*}
Since $\varphi$ is a chart, we can identify $X$ with these coordinates $c\in\mathbb R^d$ and even more also $p$ with its coordinates $\varphi(p)$. For both coordinate representations we write $X_c$ and $p_c$, respectively, for short. Similarly we can represent $X$ also with respect to a basis $B_p=\{e_1, e_2, \dots, e_d\}$ of by coefficients $c_{B_p}(X)$, that is $X = \sum_{i=1}^d c_{B_p}(X)_i e_i$.

We denote by $T\mathcal M$ the disjoint union of all tangent spaces $T \mathcal{M} = \{(p, T_p\mathcal M) \mid p \in \mathcal M\}$, called the \term{tangent bundle}.

Similar to a point in coordinates, we can also consider the curve $\gamma_c(t) \coloneqq \varphi(\gamma(t))$ in coordinates of the chart as well as its derivative $\deriv{t}{\gamma_c}$ in the same way as a tangent vector. More generally we consider \term{vector fields along a curve} $\Xi(t) \in T_{\gamma(t)}\mathcal M$, which we can again also write in coordinates denoted by $\Xi_c(t) \in \mathbb R^d$.

Additionally we assume that $\mathcal{M}$ is equipped with an affine connection, which is a method of relating tangent spaces at different points.
The relation is described by the \term{Christoffel symbol} $\Gamma_i \colon \mathbb{R}^d \times \mathbb{R}^d \times \mathbb{R}^d \to \mathbb{R}^d$~\cite[Def.~15.2]{Gallier:2020} which depends on the chart $\varphi$ as well.

A vector $Y \in T_p \mathcal{M}$ can undergo parallel transport in direction $X \in T_p \mathcal{M}$ by solving the following ODE in a chart $\varphi$ around $p$:\\
Find the curve $\gamma$ and the vector fields $\Xi, \Psi$ along $\gamma$ such that
\begin{equation}
    \label{eq:parallel_transport}
    \deriv{t}{
    \begin{pmatrix}
        \gamma_c(t)\\
        \Xi_c(t)\\
        \Psi_c(t)
    \end{pmatrix}} = \begin{pmatrix}
        \Xi_c(t) \\
        -\Gamma(\gamma_c(t), \Xi_c(t), \Psi_c(t)) \\
        -\Gamma(\gamma_c(t), \Xi_c(t), \Psi_c(t)),
    \end{pmatrix},
\end{equation}
With the initial conditions
$\gamma_c(0)=p_c$, $\Xi_c(0) = X_c$, and $\Psi_c(0) = Y_c$.

Then, the \term{parallel transport} of $Y$ in direction $X$ is then defined as the tangent vector
$\Psi(1) \in T_{\gamma(1)}$ and denoted as $\operatorname{PT}_{p,X} Y$.
Under mild conditions the computation of parallel transport does not depend on the selection of the chart. For more details see \cite[p.~105]{Lee:2018}, ~\cite[Sec.~10.3]{Boumal:2023}.

Two further concepts can be constructed based on parallel transport.
The first is the \term{exponential map} $\exp_p\colon T_p \mathcal{M} \to \mathcal{M}$ which, for a tangent vector $X$ in $T_p \mathcal{M}$, returns $\gamma(1)$ obtained by parallel transport of $X$ indirection of $X$. If this solution exists for all $X \in T_p\mathcal M$, the manifold is called \emph{complete}\cite[p.~166]{Lee:2018}. Informally, geodesics are curves on the manifold that are acceleration-free w.r.t.\ the Riemannian metric.

Locally around the origin of $T_p\mathcal M$, the exponential map is bijective. The corresponding inverse is called the \term{logarithmic map} $\log_p q$, $q\in \mathcal{M}$, that is, it finds the vector $X \in T_p \mathcal{M}$ such that $\exp_p X = q$.

Exponential and logarithmic maps are sometimes replaced by their local approximations called, respectively, retractions and inverse retraction.
A retraction is a function $\operatorname{retr}\colon T_p \mathcal{M} \to \mathcal{M}$ such that $\operatorname{retr}(0) = p$ and the differential of $\operatorname{retr}$ at the zero vector is the identity function on $T_p\mathcal{M}$.
An inverse retraction is an inverse function to a retraction whose domain is restricted to an open subset of $T_p \mathcal{M}$ containing the zero vector and on which it is bijective.

The structure provided by an affine connection is sufficient for some algorithms~\cite{PennecLorenzi:2020} but in many cases we can use additional structure, such as Lie group or a Riemannian metric.

A Lie group $\mathcal{G}$ is a manifold with a smooth group operation $\circ \colon \mathcal{G} \times \mathcal{G} \to \mathcal{G}$, where we denote the identity element by $e_{\mathcal{G}}$, and smooth inverse function $\cdot^{-1} \colon \mathcal{G} \to \mathcal{G}$ satisfying the standard group conditions:
\begin{enumerate}
    \item for all $a, b, c \in \mathcal{G}$ we have associativity, $(a \circ b) \circ c = a \circ (b \circ c)$
    \item for all $a \in \mathcal G$ we have $a \circ e_{\mathcal{G}} = e_{\mathcal{G}} \circ a = a$, and
    \item $a \circ a^{-1} = a^{-1} \circ a = e_{\mathcal G}$ holds for all $a \in \mathcal G$.
\end{enumerate}
The group structure leads to a family of affine connections, called Cartan-Schouten connections~\cite{PennecLorenzi:2020}.
Notably all such structures lead to the same exponential and logarithmic maps.

Riemannian manifolds impose additional structure by imposing an inner product $\langle \cdot, \cdot \rangle_p \colon T_p \mathcal{M} \times T_p\mathcal{M} \to \mathbb{R}$ on each tangent space that varies smoothly between points called the \term{Riemannian metric}~\cite[Def.~3.5.3]{Boumal:2023}.
On each Riemannian manifold there is a unique metric-preserving, torsion-free affine connection, called the Levi-Civita connection, see Section 15 of~\cite{Gallier:2020}.
We can also use the inner product to define a norm on tangent spaces as $\lVert X \rVert_p = \sqrt{\langle X, X \rangle_p}$ for $X \in T_p \mathcal{M}$.

Affine connections, group operations and Riemannian metrics can be broadly referred to as structure the manifold is equipped with.
It is worth noting that commonly there are multiple different ways of imposing structure required by specific algorithms.
This structure is imposed by both constraints on the configuration space $\mathcal M$ as well as by properties of $f$ from Eq.~\eqref{eq:f_system}.
The problem that is solved needs to be investigated to select the appropriate mathematical structure for the state space, that is the set of all possible values of parameters that describe the modelled system.
Group structure is typically considered to be more important than Riemannian structure for state estimation problems~\cite{SolaDerayAtchuthan:2021,BarrauBonnabel:2018}.
Notably, only Lie groups which are direct products of compact Lie groups and a Euclidean space have a Cartan-Schouten connection that is a Levi-Civita connection of some Riemannian metric.
For example, special Euclidean groups are not such direct products and thus we have to choose between a Cartan-Schouten connection or a Levi-Civita connection.

\subsection{Jacobians}

Given a function $f\colon \mathcal{M}_1 \to \mathcal{M}_2$ where $\mathcal{M}_1$ and $\mathcal{M}_2$ are manifolds, its differential at $p \in \mathcal{M}_1$ is defined as the linear map $Df(p)\colon T_{p}\mathcal{M}_1 \to T_{f(p)} \mathcal{M}_2$ such that for a chart $\varphi$ on $\mathcal{M}_2$ whose domain contains $f(p)$, and for a tangent vector $X$ at $p$ represented by a curve $\gamma$, the value of $Df(p)[X]$ is represented by the derivative of the curve $t \mapsto \varphi^{-1}(f(\gamma(t)))$.
Kalman filters involve calculation of Jacobians which are matrix representations of the linear map $Df(p)$ in a basis $B_1$ of $T_p \mathcal{M}_1$ and $B_2$ of $T_{f(p)} \mathcal{M}_2$.

In robotics literature left, right or crossed Jacobians are commonly encountered~\cite{SolaDerayAtchuthan:2021}.
They are defined for functions $f$ between two Lie groups.
They actually are standard Jacobians of the function composed with group operation as follows.
\begin{enumerate}
    \item The right Jacobian of $f$ at $p$ is the standard Jacobian of $f$ at $p$.
    \item The left Jacobian of $f$ at $p$ is the standard Jacobian of $q \mapsto f(q \circ p) \circ f(p)^{-1}$ evaluated at $q=e_{\mathcal{M}_1}$.
    \item The crossed Jacobians of $f$ at $p$ are the standard Jacobians of either $q \mapsto f(q) \circ f(p)^{-1}$ evaluated at $q=p$ or $q \mapsto f(q \circ p)$ evaluated at $q=e_{\mathcal{M}_1}$.
\end{enumerate}
Similar concepts can be developed for any parallelizable manifold.
Sometimes operations called (left- or right-)$\oplus$ and $\ominus$ are introduced to denote exponential and logarithmic map composed with group operation.
We express the theory using generic Jacobians, exponential and logarithmic maps because we do not assume the state space to be a Lie group.

\subsection{Manifold normal distribution}
\label{sec:manifold-normal}

State estimation methods use different families of distributions to express uncertainty or take a nonparametric approach like particle filters~\cite{Groves:2013}.
In this work we focus on the case where a suitable generalization of the normal distribution to manifolds is used.

Affine connection alone is not sufficient to define densities since we need a measure with respect to which densities would be defined.
However, compatibility with additional structure restricts our choice.
A suitable, unique generalization of the Lebesgue measure exists on Riemannian manifolds~\cite{Lee:2012:1,BoyaSudarshanTilma:2003}.
On Lie groups, there are unique (up to a constant scaling factor) Haar measures~\cite{Cohn:2013}.

There are multiple different generalizations of the normal distribution to manifolds~\cite{ChevallierLiLuDunson:2022,Sommer:2020}.
The standard approach in robotics is sometimes called “concentrated Gaussian”~\cite{GeGoorMahony:2024}.
These generalizations are determined by the mean $p_\mu \in \mathcal{M}$ and a positive definite covariance matrix $\Sigma \in \mathbb{R}^{d\times d}$ expressed in basis $B_{p_\mu}$.
The formula for coefficient $\Sigma^{ij}$ reads
\begin{equation}
    \Sigma^{ij} = \int_{\mathcal{M}} c_{B_{p_\mu}}(\log_p(q))_i c_{B_{p_\mu}}(\log_p(q))_j \mu(dq),
\end{equation}
where integration is performed with respect to the probability measure $\mu$ corresponding to the selected distribution~\cite{PennecLorenzi:2020}.
Its inverse is a bilinear form on $T_p\mathcal{M}$~\cite{PennecLorenzi:2020}.

The difference between different manifold normal distributions is in the formula for probability density.
In the case of “concentrated Gaussian”, the formula for the density $p_{pdf,c} \colon \mathcal{M} \to \mathbb{R}_{\geq 0}$ reads
\begin{equation}
    p_{pdf,c}(p; p_\mu, \Sigma) = A \exp\left(-\tfrac{1}{2} c_X^{\mathrm{T}} \Sigma c_X\right),
\end{equation}
where $c_X=c_{B_{p_\mu}}(\log_{p_\mu}(p))$ and $A$ is a normalization constant.
The probability density $p_{pdf,w} \colon \mathcal{M} \to \mathbb{R}_{\geq 0}$ of the exponential wrapped normal distribution under some mild conditions~\cite{ChevallierLiLuDunson:2022} reads
\begin{equation}
    \label{eq:exponential_wrapped_normal_pdf}
    p_{pdf,w}(p; p_\mu, \Sigma) = B \nu_{p_\mu}(p) \exp\left(-\tfrac{1}{2} c_X^{\mathrm{T}} \Sigma c_X\right),
\end{equation}
where $\nu_{p_\mu}\colon \mathcal{M} \to \mathbb{R}$ is volume density at point $p$ and $B$ is a different normalization constant.
The volume density term accounts for distortion caused by the exponential map~\cite{ChevallierLiLuDunson:2022}.

It is worth noting that using the exponential-wrapped normal distribution to model uncertainty has some advantages.
It allows for easy sampling, parameter estimation and density computation~\cite{ChevallierLiLuDunson:2022}.
Sampling follows by drawing a random sample $X_c$ from $N(0,\Sigma)$, taking the vector $X$ with coefficients $c_{B_p}(X)$ in a selected basis $B_p$ of $T_p \mathcal{M}$ and computing $\exp_p X$.
This value is the sample from the exponential-wrapped normal distribution.

From the sampling algorithm we can infer two things.
First, the covariance matrix depends on the selection of basis.
Second, covariance matrices at different points are not directly comparable without some identification of tangent spaces.

\subsection{Dynamical systems}

Given an initial state $p_0 \in \mathcal M$ of the system, we could solve the differential equation~\eqref{eq:f_system} using one of many existing methods.
Its choice generally depends on the exact structure of our equation.
For example, there are standard Euclidean solvers if $\mathcal{M}$ is flat, Runge-Kutte-Munte-Kaas (RKMK) methods for Lie groups~\cite{Munthe-Kaas:1998}, Crouch-Grossman methods for problems that can be expressed in the frozen coefficient formulations~\cite{CrouchGrossman:1993} or methods for differential-algebraic systems of equations which can be used if $\mathcal{M}$ is a level set of a Euclidean space, see~\cite[Section 4.5]{KunkelMehrmann:2006}.

We need to deal with an uncertainty in the initial state and random factors that can influence our system.
The first can be described using some distribution over $\mathcal{M}$ that can be evolved, while the second can be accounted for by extending the formulation to a stochastic differential equation.
We simplify this approach by applying time discretization of the system first, and then adding stochasticity.
In vast majority of applications a forward Euler discretization of the system turns out to be sufficient~\cite{Groves:2013}.
Computational effort is then spent on accurate approximation of the distribution on the space of states, in short referred to as state space.

Further discussion will be focused on discrete filtering.
This allows us to escape considering interactions between non-additive noise and curvature of the state space in the continuous setting.

Forward Euler discretization with a constant step size $\Delta t$ leads to a discretized dynamics function $\tilde{f}$ defined by
\begin{equation}
\label{eq:discretization_first_order}
    \tilde{f}(p, q, t) = \exp_p(\Delta t f(p, q, t)).
\end{equation}
For subsequent discussion we will add another parameter to $\tilde{f}$, denoted $w\in \mathbb{R}^d$, that corresponds to noise.
To summarize, discrete stochastic dynamics is thus defined by a function
\begin{equation}
\label{eq:tilde_f} 
\tilde{f}\colon \mathcal{M} \times Q \times \mathbb{R}^d \times [0, T] \to \mathcal{M}.
\end{equation}

\section{Methods}

This section presents a new unified description of discrete Kalman filters on manifolds with an affine connection.
Next, specializations to different particular variants are presented and relations between them are discussed.


Discrete Kalman filters estimate the state of a discrete stochastic system based on previous estimates, incorporating measurements.
A discrete stochastic system consists of the following ingredients:
\begin{itemize}
    \item Configuration space $\mathcal M$.
    \item Initial state $p_0 \in \mathcal M$.
    \item Stochastic dynamics function $\tilde{f}$, see Eq.~\eqref{eq:tilde_f}, with time step $\Delta t$, and final state $p_K$ for an integer $K>0$ where $T=K\Delta t$ is the total time of the simulation.
    \item Sequence of control parameters $q_0, q_1, \dots, q_{K-1} \in Q$.
    \item Sequence of process noise values $w_0, w_1, \dots, w_{K-1} \in \mathbb{R}^d$ that come from a zero-mean distribution with covariance $Q_0,\ldots Q_{K-1}$ one for each step step $n\in \{0, \dots, K-1\}$.
    \item Sequence of unobserved states $p_1, \dots, p_K$ where $p_{n+1}=\tilde{f}(p_n, q_n, w_n, n\Delta t)$ for $n = 0, 1, \dots, K-1$.
    \item Measurement function
    \begin{equation}
        \label{eq:h} 
    h\colon \mathcal{M} \times Q \times \mathbb{R}^{d_{\mathcal N}} \times [0, T] \to \mathcal{N}
    \end{equation}
    for some $d_{\mathcal N}$-dimensional manifold $\mathcal{N}$ called the measurement space.
    \item Sequence of observation noise values $v_1, v_2, \dots, v_K \in \mathbb{R}^{d_{\mathcal N}}$ that come from a zero-mean distribution with covariances $R_0,\ldots,R_{K}$, one for each step $n\in \{1, \dots, K\}$.
    \item Sequence of observations $z_n = h(p_n, q_n, v_n, n\Delta t) \in \mathcal{N}$ for $n\in\{1, \dots, K\}$.
\end{itemize}

In the filtering problem we need to estimate $p_n$ for $n=0, 1, \dots, K$ based on $\tilde{f}$, $q_n$, $h$, $z_n$ and a prior distribution for $p_0$.
Variants of this problem often add more assumptions about the distributions from which $v_n$ and $w_n$ are drawn, although their values are not known to the estimation procedure.

\subsection{Discrete Kalman filters}

Discrete Kalman filters track the evolution of a discrete stochastic system described by a state transition function $\tilde{f}$ from \eqref{eq:tilde_f} and a measurement function $h$ from \eqref{eq:h}.
The input data consists of the initial state $p_0$, its uncertainty $P_0 \in \mathbb R^{d\times d}$, a covariance matrix at $p_0$, process covariance matrices $Q_0, \dots, Q_{K-1}$, the sequence of observations $z_1, \dots, z_K$ and measurement covariance matrices $R_1, \dots, R_K$.

Kalman filters are methods that produce a sequence of state estimates $p_{0|0},\allowbreak p_{1|0},\allowbreak p_{1|1},\allowbreak \dots,\allowbreak p_{n-1|n-1},\allowbreak p_{n|n-1},\allowbreak p_{n|n},\allowbreak \dots,\allowbreak p_{K|K} \in \mathcal{M}$ and their uncertainties $P_{0|0},\allowbreak P_{1|0},\allowbreak P_{1|1},\allowbreak \dots,\allowbreak P_{n-1|n-1},\allowbreak P_{n|n-1},\allowbreak P_{n|n},\allowbreak \dots,\allowbreak P_{K|K} \in \mathbb R^{d\times d}$.
For each $n \in \{0, 1, \dots, K\}$ the point $p_{n|n} \in \mathcal{M}$ denotes the state estimate based on $p_0$, $P_0, Q_0, \dots Q_n, z_1, \dots, z_K$ and $R_1, \dots, R_n$.
The corresponding uncertainty based on the same information is denoted by $P_{n|n} \in \mathbb R^{d\times d}$, a covariance matrix at $p_{n|n}$.
Moreover, for $n \in \{1, \dots, K\}$ the point $p_{n|n-1} \in \mathcal{M}$ denotes the state estimate based on $p_0$, $P_0, Q_0, \dots Q_n, z_1, \dots, z_{n-1}$ and $R_1, \dots, R_{n-1}$.
The corresponding uncertainty based on the same information is denoted by $P_{n|n-1} \in \mathbb R^{d\times d}$, a covariance matrix at $p_{n|n-1}$.
Thus, our Bayesian estimate of the state of the discrete stochastic system at time step $n \in \{0, 1, \dots, K\}$ after incorporating observations $z_0, \dots, z_m$ for $m \in \{n-1, n\}$ is described by a chosen manifold generalization of the normal distribution with mean $p_{n|m}\in \mathcal{M}$ and covariance matrix $P_{n|m} \in \mathbb R^{d\times d}$.

The two main parts of a Kalman filter are prediction and update steps.
For $n \in \{1, 2, \dots, K\}$ the prediction steps computes a new predicted value $p_{n|n-1}$ and covariance matrix $P_{n|n-1}$, while the update steps compute values $p_{n|n}$ and covariance matrix $P_{n|n}$.
Prediction step represents evolving the initial state and uncertainty through a discretization of the dynamical system $\tilde f$ (see Eq.~\eqref{eq:tilde_f}) together with forward sensitivity analysis based on process covariance $Q_{n-1}$.


The update step entails incorporating observation $z_n \in \mathcal{N}$ with uncertainty expressed by covariance matrix $R_n$ into the state estimate.
The observation is a value of the measurement function $h$ from Eq.~\eqref{eq:h} at a particular time.
The value of observation we expect based on the state of the Kalman filter at $p_{n|n-1}\in \mathcal M$ is
\begin{equation}
o_{\mathrm{e}} = h(p_{n|n-1}, q_n, v_n, t_n) \in \mathcal N,
\end{equation}
where $v_n=0 \in \mathbb R^{d_{\mathcal N}}$ is the zero vector.
Next, we calculate measurement residual as the difference between the value of observation we expect and the actual observation:
\begin{equation}
    y_n = \log_{o_{\mathrm{e}}}(z_n).
\end{equation}
Note that $y_n$ is a tangent vector from $T_{o_{\mathrm{e}}} \mathcal{N}$.

The new state estimate is given by the maximizing the posterior probability density.
Its value can be calculated using Bayes theorem~\cite{SkoglundHendebyAxehill:2015}:
\begin{equation}
    \label{eq:Bayes_Kalman}
    p_{\mathrm{pdf}}(q|z_n) \propto p_{\mathrm{pdf}}(z_n|q) p_{\mathrm{pdf}}(q|p_{n|n-1}),
\end{equation}
where $p_{\mathrm{pdf}}$ is the corresponding conditional probability density of the selected manifold normal distribution.
For concentrated Gaussian distribution this corresponds to the following optimization problem~\cite{HuaiGao:2023}:
\begin{equation}
    \label{eq:update_kf_problem}
    p^* = \arg\min_{p\in \mathcal{M}} c_1(p)^{\mathrm{T}} P_{n|n-1}^{-1} c_1(p) + c_2(p)^{\mathrm{T}} R_{n}^{-1} c_2(p),
\end{equation}
where
\begin{equation*}
    \begin{aligned}
        c_1(p) &= c_{B_{p_{n|n-1}}}(\log_{p_{n|n-1}}p),\\
        c_2(p) &= c_{B_{o_{\mathrm{e}}}}(\log_{o_{\mathrm{e}}} h(p, q_n, 0, t_n))
    \end{aligned}
\end{equation*}
and $R_{n}$ is the covariance matrix of the observation noise.
Note that this problem depends on neither the selection of basis of either tangent space, nor the presence (or absence) of Riemannian metric or group operation on $\mathcal{M}$ or $\mathcal{N}$.

Using exponential wrapped Gaussian results in a different optimization problem.
Substituting Eq.~\eqref{eq:exponential_wrapped_normal_pdf} into Eq.~\eqref{eq:Bayes_Kalman}, taking logarithm and dropping constants gives
\begin{equation}
\begin{split}
    p^* = & \arg\min_{p\in \mathcal{M}} \left( c_1(p)^{\mathrm{T}} P_{n|n-1}^{-1} c_1(p) + c_2(p)^{\mathrm{T}} R_{n}^{-1} c_2(p) \right. \\
    & \left. -2\log(\nu_{p_{n|n-1}}(p)) - 2\log(\nu_{o_{\mathrm{e}}}(h(p, q_n, 0, t_n)) \right).
\end{split}
\end{equation}
On flat manifolds volume density $\nu_{q}(p)$ is equal to 1 for all $p,q$ and the additional term vanish.
In general, however, these terms are not negligible.
In practice measurement space is usually flat, so for simplicity we will only keep the first volume density correction term in the following text.

We can either consider this modified problem as a general nonlinear optimization problem or simplify it by replacing logarithms of volume density with its quadratic approximations around $p_{n|n-1}$:
\begin{equation}
    \begin{split}
        p^* = & \arg\min_{p\in \mathcal{M}} \big( c_1(p)^{\mathrm{T}} P_{n|n-1}^{-1} c_1(p) + c_2(p)^{\mathrm{T}} R_{n}^{-1} c_2(p) \\
        &  -2 \nabla \log(\nu_{p_{n|n-1}}({p_{n|n-1}}))^{\mathrm{T}} c_1(p) \\
        &  - c_1(p)^{\mathrm{T}} \left(\operatorname{Hess} \log(\nu_{p_{n|n-1}} (p_{n|n-1}))\right) c_1(p) \big).
    \end{split}
\end{equation}
Such simplification brings this new problem into a similar quadratic formulation as in the concentrated Gaussian case in Eq.~\eqref{eq:update_kf_problem}.

Further calculations depend on a few quantities that can be calculated in different ways that will be discussed in following sections.
The first one is the innovation covariance $S_n$ which is a tensor of rank 2 at $o_{\mathrm{e}}$ and the second one is the Kalman gain $K_n$, which is a linear operator from $T_{o_{\mathrm{e}}}\mathcal{N}$ to $T_{p_{n|n-1}} \mathcal{M}$.
The updated state estimate is now
\begin{equation}
    p_{n|n} = \exp_{p_{n|n-1}}(K_n y_n),
\end{equation}
where $K_n y_n$ is the application of $K_n$ to $y_n$.
In matrix form, we need to assure that the basis used to represent $y_n$ is the same as the basis used to represent the domain of $K_n$.
Finally, we can update the covariance matrix to the new state estimate as
\begin{equation}
P_{n|n} = \operatorname{PT}_{p_{n|n-1} \to p_{n|n}} (P_{n|n-1} - K_n S_n K_n^{\mathrm{T}}),
\end{equation}
where $\operatorname{PT}$ represents covariance transport from $p_{n|n-1}$ to $p_{n|n}$. This step is sometimes known as covariance reset, especially in the context of error-state Kalman filtering, see~\cite{GeGoorMahony:2023,MahonyGoorHamel:2022}.
Geometric formulation makes the necessity of this step on non-flat manifolds clear.
A covariance matrix can be transported by parallel transport of its eigenvectors, see Appendix~\ref{app:pt_rank_2} for a proof.

\subsubsection{Extended Kalman filter}

Mean propagation in the extended Kalman filter is simply an application of discrete dynamics to the previous state
\begin{equation}
p_{n|n-1} = \tilde{f}(p_{n-1|n-1}, q_n, w_n, t_n),
\label{eq:mean_propagation}
\end{equation}
where $p_{n|n-1}$ is the state estimate at time point $n$ given no new observation data, $q_n$ represents control, $w_n$ is noise (set to a zero vector), and $t_n$ is the time for step $n$.

Propagating uncertainty requires parametrization of the manifold using a chart  $\varphi_p \colon U \to \mathbb{R}^d$, where $U$ is a subset of $\mathcal{M}$ containing $p$. One particular choice of parametrization is related to normal coordinates. For a basis $B_p$ of $T_p \mathcal{M}$ we can define inverse of the chart as $\varphi_p^{-1}(p_c) = \exp_p \sum_{i=1}^d B_{p,i} p_{c,i}$, and thus $\varphi_{n,p}(q) = c_{B_p}(\log_p(q))$.
There are also other choices which lead to different linearization errors for uncertainty propagation. The parametrized $\tilde{f}$ reads
\begin{equation}
    \hat{f}(p_c, q_n, w_n, t_n) = \varphi_{n,2}(\tilde{f}(\varphi^{-1}_{n,1}(p_c), q_n, w_n, t_n)),
\end{equation}
where $\varphi_{n,1}$ and $\varphi_{n,2}$ may correspond to the same or different parametrization.
Note that $\hat{f} \colon \mathbb{R}^d \times Q \times \mathbb{R}^d \times [0, T] \to \mathbb{R}^d$ and the standard Jacobian with respect to the first argument can thus be computed.

Note that for some manifolds (notably the Euclidean space) it's possible to pick a single global choice of parametrization $\varphi_{n,p}$.
For compact manifolds, however, using more than one such map is necessary to avoid very high linearization errors.
The technique of moving linearization between updates is sometimes called the “error-state Kalman filter”, see~\cite{Im:2024} for a more thorough discussion.

We denote by $F_n$ the Jacobian of $\hat{f}(p_{n-1|n-1}, q_{n}, w_n, t_n)$ with respect to $p_c$ at $p_{n-1|n-1}$, that is $\varphi_{n,1}(p_{n-1|n-1}) = p_c$, and by $L_n$ the Jacobian of the same function with respect to $w$.
Now the covariance update reads
\begin{equation}
P_{n|n-1} = F_n P_{n-1|n-1} F_n^{\mathrm{T}} + L_n Q_n L_n^{\mathrm{T}},
\label{eq:covariance_propagation}
\end{equation}
where $Q_n$ is the covariance matrix of the process noise, $P_{n-1|n-1}$ is the covariance matrix at $p_{n-1|n-1}$ and $P_{n|n-1}$ is a rank 2 tensor at $p_{n|n-1}$.
This concludes the prediction step.

Update step is performed as a single step of the Gauss-Newton algorithm solving Eq.~\eqref{eq:update_kf_problem}.
Calculation of innovation covariance and Kalman gain in EKF requires the Jacobian of $h$ with respect to $p$.
Similarly to Jacobians of $\tilde{f}$, we need local parametrization $\varphi_{n,1}$ of $\mathcal{M}$ around $p_{n|n-1}$ and $\eta_{n,o_{\mathrm{e}}}$ of $\mathcal{N}$ around $o_{\mathrm{e}}$:
\begin{equation}
    \hat{h}(p_c, q_n, v_n, t_n) = \eta_{n,o_{\mathrm{e}}}(h(\varphi^{-1}_{n,1}(p_c), q_n, v_n, t_n)).
\end{equation}
Now, we denote by $H_n$ the Jacobian of $\hat{h}$ with respect to $p_c$ at $\varphi_{n,1}(p_{p|p-1})$ and by $V_n$ the Jacobian of the same function with respect to $v$.
Note that $H_n$ can be interpreted as a linear operator from $T_{p_{n|n-1}}\mathcal{M}$ to $T_{o_{\mathrm{e}}} \mathcal{N}$.
We can now define innovation covariance $S_n$ and Kalman gain $K_n$ as
\begin{align}
    S_n &= H_n P_{n|n-1} H_n^{\mathrm{T}} + W_n R_n W_n^{\mathrm{T}} \label{eq:innovation_covariance} \\
    K_n &= P_{n|n-1} H_n^{\mathrm{T}} S_n^{-1} \label{eq:kalman_gain},
\end{align}
where $R_n$ is covariance matrix of the observation noise.

A variant performing multiple steps of optimization of Eq.~\eqref{eq:update_kf_problem} is known as iterated extended Kalman filter~\cite{Im:2024,HuaiGao:2023}.
Other approaches, like Levenberg-Marquardt or quasi-Newton optimization have also been considered~\cite{SkoglundHendebyAxehill:2015,BourmaudMegretGiremusBerthoumieu:2016}.
These optimization algorithms also have manifold variants available~\cite{Boumal:2023,AdachiOkunoTakeda:2022} and can be used in this setting.

\subsubsection{Unscented Kalman filter}

Another approach to mean and covariance propagation is the unscented Kalman filter, originally proposed by~\cite{JulierUhlmann:1997} and initially extended to the manifold case by~\cite{HaubergLauzePedersen:2012}.
More recent works include~\cite{KulikovaKulikov:2022}.
A fairly thorough overview of older variants is given by~\cite{MenegazIshiharaBorgesVargas:2015}.
It replaces EKF equations with different formulas based on weighted averages while keeping the general structure of the filter unchanged.

First, we need to establish a set of sigma points around $p_{n-1|n-1}$.
Following~\cite{HaubergLauzePedersen:2012} and~\cite{WanMerwe:2000}, we define
\begin{equation}
    \begin{aligned}
        \sigma_0 &= p_{n-1|n-1} \\
        \sigma_m &= \exp_{p_{n-1|n-1}} \left( \sum_{i=1}^d f_{i,m} e_i \right)\\
        \sigma_{m+d} &= \exp_{p_{n-1|n-1}} \left( -\sum_{i=1}^d f_{i,m-d} e_i \right)
    \end{aligned}
\end{equation}
for $m \in \{1, 2, \dots, d\}$, a basis $\{e_i\}_{i=1}^d$ in which $P_{n-1|n-1}$ is expressed and a matrix $[f_{i,m}]_{d\times d} = \sqrt{(d+\lambda) P_{n-1|n-1}}$ for some scaling parameter $\lambda = \alpha^2(d+\kappa) - d$, $\alpha > 0$ and $\kappa \geq 0$.
Square root of this positive definite matrix is understood as taking its lower-triangular Cholesky factor $A$, that is $AA^{\mathrm{T}} = (d+\lambda) P_{n-1|n-1}$.

In the original unscented Kalman filter the new state $p_{n|n-1}$ is then calculated as a weighted average of points $\{\tilde{f}(\sigma_i, q_n, 0, t_n)\}_{i=0}^{2d}$.
A suitable generalization to manifolds with affine connection is given by exponential barycenter~\cite{PennecArsigny:2013}, which is a solution to the equation
\begin{equation}
    \sum_{i=0}^{2d} w_{m,i}\log_{p_{n|n-1}} \tilde{f}(\sigma_i, q_n, 0, t_n) = 0,
\end{equation}
where $w_{m,0} = \frac{\lambda}{\lambda+d}$ and $w_{m,i}=\frac{1}{2(\lambda+d)}$ for $i \in \{1, 2, \dots, 2d\}$.
Covariance matrix is estimated as
\begin{equation}
    P_{n|n-1} = \sum_{i=0}^{2d} w_{c,i} c_{B_{p_{n|n-1}}}(X_i) c_{B_{p_{n|n-1}}}(X_i)^{\mathrm{T}},
\end{equation}
where $X_i=\log_{p_{n|n-1}}(\tilde{f}(\sigma_i, q_n, 0, t_n))$, $w_{c,0} = \frac{\lambda}{\lambda+d} + (1-\alpha^2 + \beta)$ and $w_{c,i}=\frac{1}{2(\lambda+d)}$ for $i \in \{1, 2, \dots, 2d\}$.

The update step requires calculation of the expected observation value first.
We use the same sigma points as in the prediction step and a similar approach to weighted average, though this time in measurement manifold $\mathcal{N}$.
The expected observation is the solution to
\begin{equation}
    \sum_{i=0}^{2d} w_{m,i}\log_{o_{\mathrm{e}}} h(\sigma_i, q_n, 0, t_n) = 0,
\end{equation}
with the same weights as in prediction.
Innovation covariance $S_n$ is estimated by
\begin{equation}
    S_{n} = W_n R_n W_n^{\mathrm{T}} + \sum_{i=0}^{2d} w_{c,i} c_{B_{o_{\mathrm{e}}}}(Y_i) c_{B_{o_{\mathrm{e}}}}(Y_i)^{\mathrm{T}},
\end{equation}
where the first summand is the same as in EKF, $Y_i=\log_{o_{\mathrm{e}}}(h(\sigma_i, q_n, 0, t_n))$ and the weights are the same as in the prediction step.
Kalman gain is calculated as $K_n = P_{XY,n} S_n^{-1}$ where the cross-covariance matrix is estimated using
\begin{equation}
    P_{XY,n} = \sum_{i=0}^{2d} w_{c,i} c_{B_{p_{n|n-1}}}(X_i) c_{B_{o_{\mathrm{e}}}}(Y_i)^{\mathrm{T}}.
\end{equation}

\subsection{Specializations to particular variants}

Presented description of Kalman filters generalizes many existing filter types.
The standard extended Kalman filter can be obtained by setting $\mathcal{M}$ and $\mathcal{N}$ to be Euclidean spaces.
In this case we can select $\varphi_{n,1}(q) = \varphi_{n,2}(q) = q$ and, similarly, $\eta_{n,1}(q) = q$.
Euclidean error-state Kalman filter can be obtained by using $\varphi_{n,1}(q) = \varphi_{n,2}(q) = q - p_{n|n-1}$.

Multiplicative Kalman filter~\cite{LeffertsMarkleyShuster:1982} is obtained by selecting a matrix Lie group as the state space, such as the unitary quaternionic group $\mathrm{U}(1, \mathbb{H})$.

Invariant Kalman filters (IKF)~\cite{BarrauBonnabel:2017,BarrauBonnabel:2018} are error-state Kalman filters when $\mathcal{M}$ is a Lie group with a Cartan-Schouten connection.
The main insight behind IKF is that for some invariant dynamics (called group-affine systems), Jacobian calculation can be simplified.
Equivariant filters~\cite{MahonyGoorHamel:2022} extend this type of calculation simplification to a wider class of systems with symmetries in either state or control parameter.

\subsection{Noise adaptation for Kalman filters}

Importance of determination of covariance matrices $Q$ and $R$ are known for a long time~\cite{Mehra:1970, Heffes:1966}.
Common methods of noise adaptation can be classified into Bayesian, maximum likelihood, correlation and covariance matching approaches~\cite{Mehra:1972}.
While measurement noise can often be reliably determined beforehand, process noise is harder to determine.
Moreover, sometimes the process is not stationary.
Different adaptation strategies have been proposed to address these issues.
For example~\cite{AkhlaghiZhouHuang:2017} suggest an on-line scheme for updating covariance matrices, while~\cite{BerrySauer:2013} extend the classic correlation methods to extended and ensemble Kalman filters.

Noise adaptation can also be used as a part of a geometric Kalman filter.
For example, measurement and process covariance can be updated using standard covariance matching
\begin{align}
    R_{n+1} &= \alpha R_n + (1 - \alpha) (W_n^{-1}  (y_n y_n^{\mathrm{T}} + S_n)(W_n^{\mathrm{T}})^{-1} - R_n), \\
    Q_{n+1} &= \alpha Q_n + (1 - \alpha) (L_n^{-1} K_n y_n  (L_n^{-1} K_n y_n)^{\mathrm{T}}),
\end{align}
for some smoothing parameter $\alpha \in [0, 1]$.

\section{Results}

The proposed approach has been validated on two examples implemented in the \verb|GeometricKalman.jl| library.
It contains generic implementation of Kalman filters, applicable to configuration spaces that are Lie groups with invariant connections, Riemannian manifolds or more general manifolds equipped with an affine connection.
In this Section we present two examples, one for a Lie group and one for a Riemannian manifold.
The three Kalman filters used in both examples illustrate the versatility of both the theoretical framework and the implementation.

The first example is a simple car model on a plane based on an example from~\cite{BarrauBonnabel:2017}.
Its state $(p, X)$ is an element of the special Euclidean group $\mathrm{SE}(2)$ with the torsion-free Cartan-Schouten connection.
The discrete dynamics are defined by $\tilde{f}_{\text{car}} \colon \mathrm{SE}(2) \times \mathbb{R} \times \mathbb{R}^3 \times \mathbb{R} \to \mathrm{SE}(2)$ such that
\begin{equation*}
    \begin{split}
    &\tilde{f}_{\text{car}}((p, X), q, w, t) = \\
    & \left(\exp_p (\Delta t X v + \sqrt{\Delta t}[w_2, w_3]), \begin{bmatrix}
        0 & -q-w_1\\
        q+w_1 & 0
    \end{bmatrix}\right),
    \end{split}
\end{equation*}
where $(p, X) \in \mathrm{SE}(2)$ is the state of the system at time $t$, $q$ is the system control parameter given by $q = q(t) = \sin(t / 2)$, $w \in \mathbb{R}^3$ is the noise at the current time step.
The measurement function $h_c\colon \mathrm{SE}(2) \times \mathbb{R} \times \mathbb{R}^2 \times \mathbb{R} \to \mathbb{R}^2$ is given by
\begin{equation}
    \label{eq:car_h}
    h_c((p, X), q, w, t) = p + w,
\end{equation}
where $w$ is the measurement noise.

A simulation was performed with time step 0.01 for 200 steps, with initial point being the identity element of $\mathrm{SE}(2)$, with process noise covariance $\operatorname{diag}(10, 100, 100)$ and measurement noise covariance $\operatorname{diag}(0.01, 0.01)$, where $\operatorname{diag}(a_1, \dots, a_n)$ is the diagonal matrix with elements $a_1, \dots, a_n$ on the diagonal.

The initial position uncertainty for Kalman filters was set to $\operatorname{diag}(0.1, 0.1, 0.1)$. Process covariance for Kalman filters was set to $\operatorname{diag}(2, 2, 2)$ and measurement covariance to $\operatorname{diag}(0.01, 0.01)$.
The extended Kalman filter for manifolds with an affine connection was compared to two other filters.
The first one is EKF with measurement and process covariance matrix adaptation according to~\cite{AkhlaghiZhouHuang:2017}, extended to support non-additive process noise, with forgetting factor $\alpha=0.99$.
The second one is the unscented Kalman filter with Wan-Merwe sigma points.


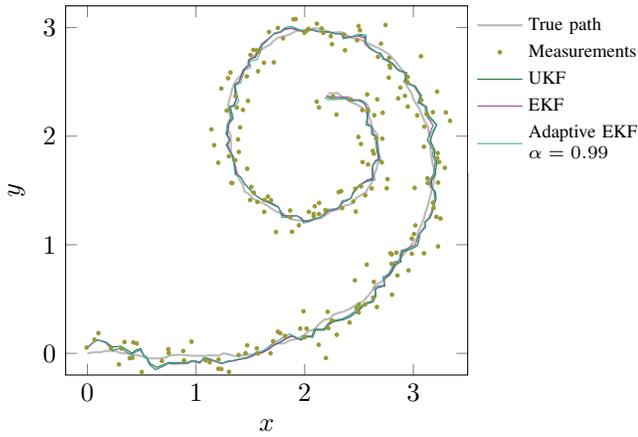
\begin{figure}[tb]  
  \centering
  \begin{tikzpicture}

    \pgfplotstableread[col sep=comma]{data/car2d_trajectory.csv}\trajdata
    \pgfplotstableread[col sep=comma]{data/car2d_measurements.csv}\measdata

    \begin{axis}[
      width=.85\columnwidth,
      xlabel={$x$},
      ylabel={$y$},
      xmin=-0.2, xmax=3.5,
      ymin=-0.2, ymax=3.2,
      xtick={0,1,...,3},
      ytick={0,1,...,3},
      axis equal image,           
      legend cell align=left,
      legend style={font=\scriptsize, at={(1,1)}, draw=none, anchor=north west, cells={align=left}},
    ]
      \addplot[TolVibrantGray, thick]
        table[x index=1, y index=2] {\trajdata};
      \addlegendentry{True path}
      \addplot[
        only marks,
        mark=*,
        mark size=0.66pt,
        TolMutedOlive,
      ]
        table[x index=1, y index=2] {\measdata};
      \addlegendentry{Measurements}
      \addplot[TolMutedGreen]
        table[x index=5, y index=6] {\trajdata};
      \addlegendentry{UKF}
      \addplot[TolMutedPurple]
        table[x index=3, y index=4] {\trajdata};
      \addlegendentry{EKF}
      \addplot[TolMutedTeal]
        table[x index=7, y index=8] {\trajdata};
      \addlegendentry{Adaptive EKF\\$\alpha=0.99$}
    \end{axis}
  \end{tikzpicture}

  \caption{Car tracking example: ground truth, raw sensor measurements, and state estimation trajectories from EKF, UKF, and an adaptive EKF ($\alpha=0.99$).}
  \label{fig:car_example}
\end{figure}

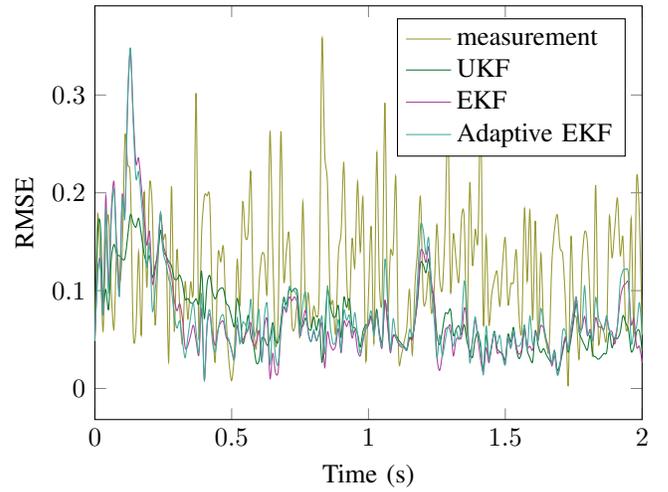
\begin{figure}[htbp]
  \centering
  \begin{tikzpicture}
    \begin{axis}[
      width=\columnwidth,            
      height=0.8\columnwidth,       
      xmin=0, xmax=2,
      ylabel={RMSE},
      xlabel={Time (s)},
      legend cell align=left,
      legend pos=north east,         
    ]

      \addplot[smooth, TolMutedOlive] table[
        col sep=comma,
        x=times,
        y=measurement_errors
      ] {data/car2d_errors.csv};
      \addlegendentry{measurement}

      \addplot[smooth, TolMutedGreen] table[
        col sep=comma,
        x=times,
        y=error_UKF
      ] {data/car2d_errors.csv};
      \addlegendentry{UKF}

      \addplot[smooth, TolMutedPurple] table[
        col sep=comma,
        x=times,
        y=error_EKF
      ] {data/car2d_errors.csv};
      \addlegendentry{EKF}

      \addplot[smooth, TolMutedTeal] table[
        col sep=comma,
        x=times,
        y={error_EKF adaptive α=0.99}
      ] {data/car2d_errors.csv};
      \addlegendentry{Adaptive EKF}

    \end{axis}
  \end{tikzpicture}
  \caption{Estimation errors for selected Kalman filter variants. The forgetting parameter of the adaptive EKF filter was set to $\alpha=0.99$.}
  \label{fig:car_example_errors}
\end{figure}

Figure~\ref{fig:car_example} shows the state tracked by different filters.
Corresponding root mean squared errors are displayed in Figure~\ref{fig:car_example_errors}.
The mean squared error of reconstructed trajectories is substantially lower than corresponding measurement error.

The second example consists of an object moving on the surface of a two-dimensional sphere.
State of the system is described by an element $(p, X)$ of the tangent bundle $T S^2$.
It can be equipped with the Sasaki Riemannian metric~\cite{Sasaki:1958}, although for computational efficiency a first order retraction and the corresponding inverse retraction were used in the experiments.

The discrete dynamics are defined by $\tilde{f}_{\text{sphere}} \colon T S^2 \times \mathbb{R} \times \mathbb{R}^4 \times \mathbb{R} \to T S^2$ such that
\begin{equation}
    \label{eq:sphere_f_tilde}
    \begin{split}
    &\tilde{f}_{\text{sphere}}((p, X), q, w, t) = \\
    & \left(v X + c_{B_p}^{-1}([w_3, w_4]), p \times [0, qv + w_1, w_2]\right),
    \end{split}
\end{equation}
where $\times$ denotes the cross product, $(p, X)$ is the state of the system at time $t$, $q$ is the system control parameter given by $\sin(t / 2)$ and $w \in \mathbb{R}^4$ is the noise at the current time step.
The measurement function $h_s\colon T S^2 \times \mathbb{R} \times \mathbb{R}^2 \times \mathbb{R} \to S^2$ is given by
\begin{equation}
    \label{eq:sphere_h}
    h_s((p, X), q, w, t) = \exp_p \left(c_{B_p}^{-1}(w) \right),
\end{equation}
where $w$ is the measurement noise, $B_p$ is a selected orthonormal basis of $T_p S^2$.

A simulation was performed with time step 0.01 for 200 steps, with initial point $p_0 = ([1, 0, 0], [0, 1, 0])$, with process noise covariance $\operatorname{diag}(10, 10, 1, 1)$ and measurement noise covariance $\operatorname{diag}(0.01, 0.01)$.

The initial position uncertainty for Kalman filters was set to $\operatorname{diag}(0.1, 0.1, 0.1, 0.1)$. Process covariance for Kalman filters was set to $\operatorname{diag}(0.1, 0.1, 0.01, 0.01)$ and measurement covariance to $\operatorname{diag}(0.01, 0.01)$.
Extended Kalman filter for manifolds with an affine connection was compared to two other filters.
The first one is EKF with measurement covariance matrix adaptation with forgetting factor $\alpha=0.99$.
The second one is the unscented Kalman filter with Wan-Merwe sigma points.
\comment[id=RB]{Also here, I am mising the punch line of the example; sorry by changing colors the ais in the 3D plot are not yet so nicem we could also move that one to an asymptote render (see Manopt); or Makie – but maybe lets first fix the story before I spent 2 days trying to get a nice 3D plot and then we do not use it.}
\begin{figure}[tb]
  \centering
  \begin{tikzpicture}

    \pgfplotstableread[col sep=comma]{data/sphere_trajectory.csv}\trajdata
    \pgfplotstableread[col sep=comma]{data/sphere_measurements.csv}\measdata

    \begin{axis}[
      view={60}{25},                
      width=0.85\columnwidth,
      scale only axis=true,
      enlarge x limits=false,   
      enlarge y limits=false,
      enlarge z limits=false,
      axis equal image,                   
      xlabel={$x$},
      ylabel={$y$},
      zlabel={$z$},
      legend style={
        font=\scriptsize,
        draw=none,
        at={(1.00,1.00)}, anchor=north west,
        cells={align=left, anchor=west}
      },
      tick label style={font=\scriptsize},
      label style={font=\scriptsize},
      xmin=-1.01, xmax=1.01,
      ymin=-1.01, ymax=1.01,
      zmin=-1.01, zmax=1.01,
      clip=false,
    ]

      \addplot3[
        domain=0:360,
        y domain=-90:90,
        samples=60,
        samples y=30,
        surf,
        shader=flat,
        draw=gray!50,
        fill=gray!10,
        z buffer=sort,
        opacity=0.1,
        forget plot               
      ]
      ({cos(y)*cos(x)}, {cos(y)*sin(x)}, {sin(y)}); 

      \addplot3[TolVibrantGray, thick]
        table[x index=1, y index=2, z index=3] {\trajdata};
      \addlegendentry{True path}

      \addplot3[
        only marks,
        mark=*,
        mark size=0.66pt,
        TolMutedOlive,
      ]
        table[x index=1, y index=2, z index=3] {\measdata};
      \addlegendentry{Measurements}

      \addplot3[TolMutedGreen]
        table[x index=7, y index=8, z index=9] {\trajdata};
      \addlegendentry{UKF}
      \addplot3[TolMutedPurple]
        table[x index=4, y index=5, z index=6] {\trajdata};
      \addlegendentry{EKF}
      \addplot3[TolMutedTeal]
        table[x index=10, y index=11, z index=12] {\trajdata};
      \addlegendentry{Adaptive EKF\\$\alpha=0.99$}

    \end{axis}
  \end{tikzpicture}

  \caption{3-D tracking on the unit sphere: ground-truth trajectory, raw sensor measurements, and state estimation results from EKF, UKF, and an adaptive EKF ($\alpha=0.99$).}
  \label{fig:sphere_example}
\end{figure}
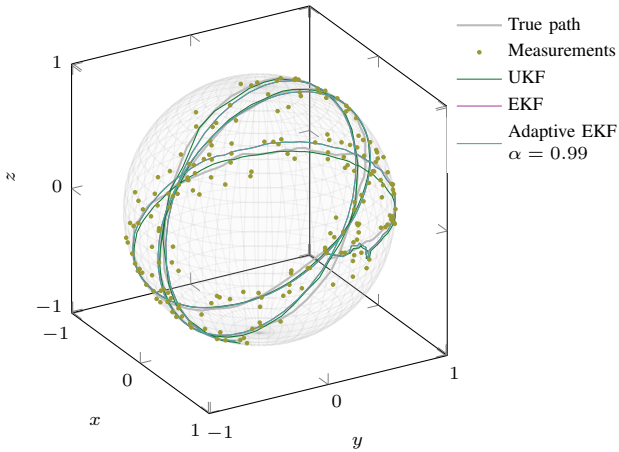

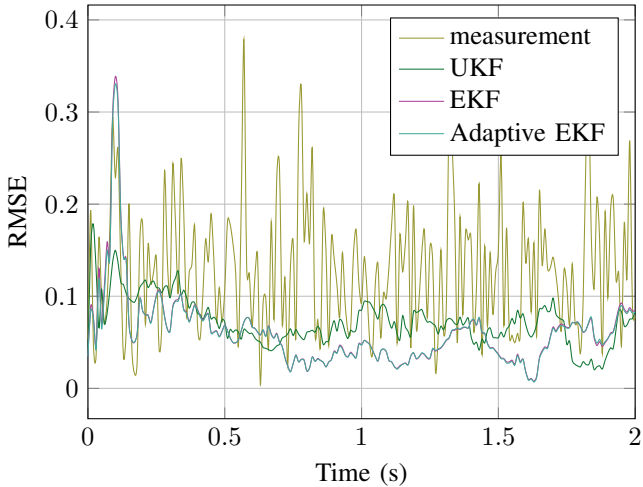
\begin{figure}[htbp]
  \centering
  \begin{tikzpicture}
    \begin{axis}[
      width=\columnwidth,            
      height=0.8\columnwidth,       
      xmin=0, xmax=2,
      ylabel={RMSE},
      xlabel={Time (s)},
      grid=major,
      legend cell align=left,
      legend pos=north east,         
    ]

      \addplot[smooth, TolMutedOlive] table[
        col sep=comma,
        x=times,
        y=measurement_errors
      ] {data/sphere_errors.csv};
      \addlegendentry{measurement}

      \addplot[smooth, TolMutedGreen] table[
        col sep=comma,
        x=times,
        y=error_UKF
      ] {data/sphere_errors.csv};
      \addlegendentry{UKF}

      \addplot[smooth, TolMutedPurple] table[
        col sep=comma,
        x=times,
        y=error_EKF
      ] {data/sphere_errors.csv};
      \addlegendentry{EKF}

      \addplot[smooth, TolMutedTeal] table[
        col sep=comma,
        x=times,
        y={error_EKF adaptive M α=0.99}
      ] {data/sphere_errors.csv};
      \addlegendentry{Adaptive EKF}

    \end{axis}
  \end{tikzpicture}
  \caption{Estimation errors for selected Kalman filter variants. The forgetting parameter of the adaptive EKF filter was set to $\alpha=0.99$.}
  \label{fig:sphere_example_errors}
\end{figure}

Figure~\ref{fig:sphere_example} shows the state tracked by different filters.
Corresponding root mean squared errors are displayed in Figure~\ref{fig:sphere_example_errors}.
The mean squared error of reconstructed trajectories is substantially lower than corresponding measurement error.

In the selected cases the results are identical to existing Kalman filters, namely Lie group filters for the first example and Riemannian filters for the second one.
The primary benefit of our new approach is to demonstrate that the same filtering logic can be applied to both scenarios, and even more complicated ones involving both Lie groups and more general manifolds.
It is sufficient to supply the right retraction, inverse retraction, vector transport, and a basis for each tangent space.

\section{Conclusions}

We have proposed a novel, generalized description of Kalman filters on manifolds equipped with an affine connection.
The results demonstrate that this approach can be successfully applied to both Lie groups and Riemannian manifolds.
This flexibility of Kalman filters can be exploited to design new variants tailored to different applications such as human motion tracking or tracking the state of under-actuated robots.
The library \verb|GeometricKalman.jl|, published together with this work, can serve as a tool for developing such solutions.

Future work includes developing second order Kalman filters~\cite{Jazwinski:1970} on manifolds using jets~\cite{KolarSlovakMichor:1993}.
Continuous variants of Kalman filters would allow for application of higher order manifold integration schemes for the prediction step.
The choice of the generalization of the normal distribution to manifolds also needs requires further exploration.

\begin{appendices}

\section{Examples of manifolds and group actions}

Tables~\ref{tab:manifolds-basics} and~\ref{tab:manifolds-functions} collect basic properties of manifolds discussed in this work.
A torsion-free connection was selected for all Lie groups.
They can also be equipped with a different Cartan-Schouten connection that has the same formulas for the exponential and logarithmic maps but differ by parallel transport.
In particular, any group can be equipped with a flat connection in which $\operatorname{PT}_{p,Y}X = X$.

Tangent bundle is not included there because it requires a relatively lengthy introduction that is beyond the scope of this manuscript.
A computationally oriented description can be found in~\cite{MuralidharanFletcher:2012}.

\begin{table*}
\caption{Selected manifolds and their properties. Tangent vectors to Lie groups are represented as elements of the Lie algebra using the differential of left translation ($\lambda_p (q) = p \circ q$) as the isomorphism between tangent spaces and the Lie algebra.}
\centering
\label{tab:manifolds-basics}
\begin{tabular}{p{90pt}p{110pt}p{90pt}p{30pt}p{60pt}p{50pt}}
\toprule
\multicolumn{3}{c}{Manifold} & dimension & \multicolumn{2}{c}{Basic operations} \\
\cmidrule(r){1-3}\cmidrule(lr){5-6}
Name & $p\in \mathcal M$ & $X \in T_p\mathcal M$ & & $\langle X, Y\rangle$ & $p\circ q$\\
\midrule
Euclidean &
$p\in\mathbb{R}^n$ & $X\in\mathbb{R}^n$ & $n$ & $X^{\mathrm{T}}Y$ & $p + q$
\\[.25\baselineskip]
Sphere &
$p\in S^n \subset \mathbb{R}^{n+1}$, $\lVert p \rVert = 1$
& $X \in \mathbb{R}^{n+1}$, $p^{\mathrm{T}}X = 0$
& $n$ & $X^{\mathrm{T}}Y$ & \emph{(not a group)}
\\[.25\baselineskip]
Rotations &
$p\in \mathrm{SO}(n) \subset \mathbb{R}^{n\times n}$, $p^{\mathrm{T}}p = \mathrm{I}_n$ & %
$X \in \mathbb{R}^{n\times n}, X^{\mathrm{T}} = -X$ & %
$\frac{n(n-1)}{2}$ &
$\operatorname{tr}(X^\mathrm{T}Y)$ & $p q$
\\[.25\baselineskip]
Unit quaternions &
$p\in \mathrm{U}(1, \mathbb{H}) \subset \mathbb H$, $\lVert p \rVert = 1$ &
$X \in \mathbb{H}$, $\mathrm{Re}(X) = 0 $ &
$3$ &
$\displaystyle\sum_{i=1}^3 X_i Y_i$ & $p q$
\\[1.25\baselineskip]
Special Euclidean group&
$(t,R) \in \mathrm{SE}(n) = \mathbb{R}^n \rtimes \mathrm{SO}(n)$ &
$(v,A) \in \mathbb R^n \times T_I\mathrm{SO}(n)$ &
$\frac{n(n+1)}{2}$ &
$v^{\mathrm{T}}w + \operatorname{tr}(A^\mathrm{T}B)$ &
$(t + R u, R S)$
\\[.25\baselineskip]
\bottomrule
\end{tabular}
\label{tab1}
\end{table*}

\begin{table*}
\caption{Selected manifolds and their exponential and logarithmic maps as well as their torsion-free parallel transport.
The formulas for special Euclidean group $\mathrm{SE}(n)$ use the embedding in the general linear group $\mathrm{GL}(n+1)$, where a point $(t, R) \in \mathrm{SE}(n)$ is represented by a matrix $\left(\begin{smallmatrix}
  R & t \\
  0 & 1
\end{smallmatrix}\right)$}
\centering
\label{tab:manifolds-functions}
\setlength{\tabcolsep}{3pt}
\begin{tabular}{p{90pt}p{120pt}p{120pt}p{120pt}}
\toprule
Manifold & Exponential map $\exp_pX$ & Logarithmic map $\log_pq$ & Parallel transport $\operatorname{PT}_{p,Y}X$
\\\midrule
Euclidean & $p + X$ & $q - p$ & $X$
\\[.25\baselineskip]
Sphere &
$\cos(\lVert X \rVert_p)p + \sin(\lVert X \rVert_p)\frac{X}{\lVert X \rVert_p}$ &
$\arccos(p^{\mathrm{T}}q) \frac{q-(p^{\mathrm{T}}q) p}{\lVert q-(p^{\mathrm{T}}q) p \rVert_2}$ &
$X - \frac{\langle Y,X \rangle_p}{\lVert Y \rVert} \bigl(Y + \log_{\exp_p X} p \bigr)$
\\[.25\baselineskip]
Rotations & %
$p \exp(X)$ &
$\log(p^{\mathrm{T}}q)$ &
$\exp(-Y/2) X \exp(Y/2)$
\\[.25\baselineskip]
Unit quaternions &
$p \exp(X)$ &
$\log(p^{-1} q)$ &
$\exp(-Y/2) X \exp(Y/2)$
\\[.5\baselineskip]
Special Euclidean group &
$p \exp(X)$ &
$\log(p^{-1}q)$ &
$\exp(-Y/2) X \exp(Y/2)$
\\
\bottomrule
\end{tabular}
\label{tab2}
\end{table*}

\section{Parallel transport of covariance matrices}

\label{app:pt_rank_2}

We use Einstein summation notation in this section.
A covariance matrix is a symmetric contravariant tensor of rank 2~\cite{PennecLorenzi:2020}.
Consider a covariance matrix $A^{ij}$ at point $p$ in some chart $\varphi$ of a $d$-dimensional manifold $\mathcal{M}$.
The ODE for parallel transport of $A^{ij}$ in direction $X^k$ along the geodesic, following Section 7.2.4 of~\cite{Nakahara:2003}, reads
\begin{equation}
    \label{eq:pt_tensor}
    \frac{d}{dt}\begin{bmatrix}
        p_c(t)\\
        X_c(t)\\
        \hat{A}^{ij}(t)
    \end{bmatrix} = \begin{bmatrix}
        X_c(t) \\
        -\Gamma(p_c(t), X_c(t), X_c(t)) \\
        -\Gamma^i_{kl} X^k(t) \hat{A}^{lj}(t) - \Gamma^{j}_{kl} X^k(t) \hat{A}^{il}(t)
    \end{bmatrix},
\end{equation}
where $\hat{A}$ is a function $[0, 1] \to \mathbb{R}^d \times \mathbb{R}^d$ such that $\hat{A}(0)^{ij}=A^{ij}$ and $\hat{A}(1)$ is the transported tensor at point with coordinates $p_c(1)$ in chart $\varphi$.
The dependence of $\Gamma^i_{kl}$ on $p_c(t)$ and chart $\varphi$ is implicit for compactness.
In this notation the vector $X_c(t)$ has coordinates $X^k(t)$ and
\begin{equation}
    \Gamma(p_c(t), X_c(t), X_c(t)) = \Gamma^i_{kl} X^k(t) X^{l}(t).
\end{equation}
Other initial conditions are analogous to the parallel transport of vectors.
We show that this operation can be expressed using parallel transport of eigenvectors of $A^{ij}$.

A straightforward computation shows that if $A^{ij}$ is symmetric, then $\hat{A}^{ij}(t)$ is symmetric for all $t\in [0, 1]$.
Since $\hat{A}^{ij}(t)$ is symmetric, all of its eigenvalues are real and thus it can be decomposed as
\begin{equation}
    \hat{A}^{ij}(t)=\sum_{u=1}^{d}\lambda^u(t) Y^i_u(t) Y^j_u(t),
\end{equation}
where $Y_u(t)$ is the $u$th eigenvector and $\lambda_u(t)$ is the corresponding eigenvalue.
Parallel transport of eigenvectors in direction $X^k$ is computed by integrating the equation
\begin{equation}
    \frac{d}{dt}Y^i_u(t) = -\Gamma^i_{kl}X^k(t) Y^l_u(t),
\end{equation}
which is equivalent to the bottom part of Eq.~\eqref{eq:parallel_transport}.

If we compute the derivative of eigendecomposition of $\hat{A}^{ij}$ we get using the product rule
\begin{equation}
    \begin{split}
        \frac{d}{dt} \hat{A}^{ij}(t) &= \frac{d}{dt} \left(\lambda^u(t) Y^i_u(t) Y^j_u(t)\right) \\
        &= \left(\frac{d}{dt} \lambda^u(t) \right) Y^i_u(t) Y^j_u(t)\\
        &\qquad + \lambda^u(t) \left(\frac{d}{dt} Y^i_u(t) \right) Y^j_u(t)\\
        &\qquad + \lambda^u(t) Y^i_u(t) \left(\frac{d}{dt} Y^j_u(t)\right).
    \end{split}
\end{equation}
Substituting expression for parallel transport of the eigenvectors we further get
\begin{equation}
    \begin{split}
        \frac{d}{dt} \hat{A}^{ij}(t) &= \left(\frac{d}{dt} \lambda^u(t) \right) Y^i_u(t) Y^j_u(t)\\
        &\qquad + \lambda^u(t) \left(-\Gamma^i_{kl}X^k(t) Y^l_u(t) \right) Y^j_u(t) \\
        &\qquad + \lambda^u(t) Y^i_u(t) \left(-\Gamma^j_{kl}X^k(t) Y^l_u(t)\right).
    \end{split}
\end{equation}
Rearranging the terms on the right hand side gives
\begin{equation}
    \begin{split}
       \frac{d}{dt} \hat{A}^{ij}(t) &= \left(\frac{d}{dt} \lambda^u(t) \right) Y^i_u(t) Y^j_u(t) \\
        &\qquad -\Gamma^i_{kl}X^k(t) \lambda^u(t) Y^l_u(t) Y^j_u(t) \\
        &\qquad -\Gamma^j_{kl}X^k(t) \lambda^u(t) Y^i_u(t) Y^l_u(t).
    \end{split}
\end{equation}
This is equal to the parallel transport as defined by Eq.~\eqref{eq:pt_tensor} if and only if $\frac{d}{dt} \lambda^u(t) = 0$, thus concluding that parallel transport of the original tensor is equivalent to parallel transport of its eigenvectors while keeping corresponding eigenvalues constant.

\end{appendices}